\documentclass[reqno]{amsart}

\usepackage{amsmath}
\usepackage{amssymb}
\usepackage{amsfonts}
\usepackage{amscd}
\usepackage{mathrsfs}
\usepackage{latexsym}
\usepackage{epsfig}
\usepackage{subfigure}

\newcommand{\R}{{\mathbb R}}
\newcommand{\C}{{\mathbb C}}
\newcommand{\p}{{\mathbb P}}

\newcommand{\Z}{{\mathbb Z}}
\newcommand{\N}{{\mathbb N}}
\newcommand{\NI}{{\noindent}}
\newcommand{\QED}{\hfill$\Box$\medskip}

\newtheorem{theorem}{Theorem}[subsection]
\newtheorem{cor}[theorem]{Corollary}
\newtheorem{defn}[theorem]{Definition}
\newtheorem*{definition}{Definition}
\newtheorem{thm}[theorem]{Theorem}

\newtheorem{lemma}[theorem]{Lemma}
\newtheorem{prop}[theorem]{Proposition}

\begin{document}

\title[Compact Complex Surfaces and cscK Metrics]{Compact Complex Surfaces and Constant Scalar Curvature K\"ahler Metrics}
\author{Yujen Shu}
\bibliographystyle{plain}
\thanks{Key words: deformation equivalence, constant scalar curvature K\"ahler
metrics.} \thanks{2000 Mathematics Subject Classification: 53C25,
53C55}

\begin{abstract}
In this article, I prove the following statement: Every compact
complex surface with even first Betti number is deformation
equivalent to one which admits an extremal K\"ahler metric. In fact,
this extremal K\"ahler metric can even be taken to have constant
scalar curvature in all but two cases: the deformation equivalence
classes of the blow-up of $\p_2$ at one or two points. The explicit
construction of compact complex surfaces with constant scalar
curvature K\"ahler metrics in different deformation equivalence
classes is given. The main tool repeatedly applied here is the
gluing theorem of C. Arezzo and F. Pacard which states that the
blow-up/resolution of a compact manifold/orbifold of discrete type,
which admits cscK metrics, still admits cscK metrics.

\end{abstract}

\maketitle

\section{Introduction}
Let $(M, J)$ be a compact complex K\"ahler manifold, and $\gamma$ be
a K\"ahler class of $(M, J)$. Calabi \cite{caI} considers the
functional $\Phi(g):=\int s_g^2 dv_g$ defined on the set of all
K\"ahler metrics $g$ in the class $\gamma$, where $s_g$ denotes the
scalar curvature associated to the metric $g$. A K\"ahler metric is
called extremal if it is a critical point of $\Phi$. It has been
shown that $g$ is extremal if and only if the gradient of the scalar
curvature $s_g$ is a real holomorphic vector field. In particular,
$g$ is extremal if $s_g$ is constant. The famous Aubin-Yau theorem
\cite{au, yau} asserts that every compact complex manifold $X$ with
negative first Chern class admits a K\"{a}hler-Einstein metric in
the canonical class $-c_1(X)$, which has a negative constant scalar
curvature. Followed by applying a deformation argument due to LeBrun
and Simanca \cite{ls}, we know there exists a constant scalar
curvature K\"ahler (cscK) metric in every class near $-c_1(X)$.

There are classical obstructions to the existence of cscK metrics
related to the automorphism group $Aut(M,J)$ of $(M,J)$. The
Matsushima-Lichnerowicz theorem \cite{li, ma} states that the
identity component $Aut_0(M,J)$ of the automorphism group must be
reductive if a cscK metric exists. Later, Futaki \cite{fu} shows
that the scalar curvature of an extremal K\"ahler metric is constant
if and only if its Futaki invariant vanishes identically. Recently
Donaldson \cite{d}, Chen-Tian \cite{ch} and Mabuchi \cite{mab, mab2}
have made substantial progress on relating the existence and
uniqueness of extremal K\"ahler metrics in Hodge K\"ahler classes to
the K-stability of polarized projective varieties. In particular, it
has been shown that K-stability is a necessary condition for the
existence of cscK metrics for a polarized projective variety. The
main result of this article is the following:

\smallskip

\noindent {\bf Main Theorem} Every compact complex surface with even
first Betti number is deformation equivalent to one which admits an
extremal K\"ahler metric. In fact, this extremal K\"ahler metric can
even be taken to have constant scalar curvature in all but two
cases: the deformation equivalence classes of the blow-up of $\p_2$
at one or two points.
\smallskip

The hypothesis that $b_1$ is even is equivalent to requiring that
$(M, J)$ admits a K\"ahler metric \cite{si, t}. One of the main
tools in our proof is the major breakthrough by Arezzo and Pacard
\cite{ap}: Let $M$ be a complex manifold/orbifold of discrete type,
which admits cscK metrics. Then the blow-up/resolution of $M$ admits
cscK metrics.

To prove the main theorem, we proceed along the Kodaira dimension
$\kappa$ of $(M,J)$. The main difficulty lies in the case of
$\kappa=1$. First of all, we deal with the case of principal
$E$-bundles, where $E$ is a smooth elliptic curve. We could show
every principal $E$-bundle with even Betti number is covered by
$\mathbb{H}\times \C$, and it inherits a cscK metric from the
product metric. Analyzing the geometric structure of the surface and
the local structure near a multiple fibre, we can generalize this
result to the case of all elliptic surfaces whose each fibre has
smooth reduction. To construct an elliptic surface with cscK metrics
and positive Euler number in each deformation class, we start with
the elliptic surface $S$, which is obtained  by applying logarithmic
transform along some fibres on the trivial elliptic bundle. Under
certain choices of the base curve and the fibres on which
logarithmic transform is applied, there exists an holomorphic
involution $\iota: S\to S$. Although the quotient of $S$ by the
action of involution $\iota$ is singular, we can show that its
resolution is smooth and carries cscK metrics  by an application of
Arezzo-Pacard's result \cite{ap}. Finally,  to show elliptic
surfaces with positive Euler number of each deformation class can be
constructed in this way, we use the deformation theory of elliptic
surfaces \cite{fm}: in the case of positive Euler number, the
deformation class of an elliptic surface is determined by the
diffeomorphism type of the base orbifold and the Euler number
$\chi$.

In section 4, we explain why the complex surfaces of other Kodaira
dimensions are deformation equivalent to complex surfaces with cscK
metrics. The case of $\kappa=2$ is done as an application in
Arezzo-Pacard's paper \cite{ap} by using the Aubin-Yau theorem
\cite{au, yau} and the fact that negative first Chern class implies
the automorphism group is discrete. Similarly, the case of
$\kappa=0$ is a result of Yau's theorem \cite{yau} and the fact that
all holomorphic vector fields on a minimal complex surface of
$\kappa=0$ are parallel. The case of ruled surfaces,
$\kappa=-\infty$, is dealt with by using the result that if a rank 2
vector bundle $E$ is poly-stable, then the associated ruled surface
$\p(E)$ admits cscK metrics \cite{at}. In the end, we also show that
$\p_2 \# k \overline{\p}_2, k=1,2$, is not deformation equivalent to
a complex surface with cscK metrics by showing that the Lie algebra
of holomorphic vector fields of every compact complex surface in the
deformation class is not reductive.

Although there exists no cscK metrics on any K\"ahler class of the
blow-up of $\p_2$ at one or two points, Calabi \cite{caI}, Arezzo,
Parcard, and Singer show \cite{ap3} that there do exist extremal
K\"ahler metrics on them, and the main theorem follows.
\smallskip

{\bf Acknowledgments:} I would like to express my gratitude to
Claude LeBrun for suggesting this problem to me and for his great
guidance and continuous encouragement. And I would like to thank Li
Li for interesting conversations on these topics and the referee for
several useful comments on this manuscript.

\section{Elliptic surfaces}
In this section, we give a general description of elliptic surfaces,
and discuss some special examples: principal $E$-bundles and the
ones without singular fibres. Then we introduce one important
operation: logarithmic transform which we will use later to
investigate the geometric structures of elliptic surfaces.

\subsection{Properties of elliptic fibrations}
\begin{defn} A compact complex surface $S$ is called an elliptic
surface if there is a surjective holomorphic map $\pi: S\rightarrow
B$ over a smooth algebraic curve $B$, whose general fibres are
smooth elliptic curves. The map $\pi$ is called an elliptic
fibration.
\end{defn}

\NI Throughout this section, we only consider the case that the
elliptic surface $S$ is minimal. That is, it contains no rational
curves of self-intersection number -1. An elliptic surface may have
singular, reducible, multiple fibres. Kodaira \cite{koII} has
classified all possible fibres of an elliptic surface, which are
denoted by $_mI_k$, $II$,$ III$, $IV, I_k^*$, $II^*, III^*$, and
$IV^*$, where $k\in \N\cup \{0\}$. The set of multiple fibres plays
a central role in our discussion of elliptic surfaces. And we will
recall the definition here.
\begin{defn}
Let $\pi: S\rightarrow B$ be an elliptic fibration and $p\in B$. We
call $\pi^{-1}(p)$ a multiple fibre if there exists an integer
$m\geq 1$ and a reduced divisor $D$ in $S$ such that as a divisor
$\pi^{-1}(p)=mD$. The largest $m$ is called the multiplicity of the
fibre.
\end{defn}

In fact, there are only two types of multiple fibres in an elliptic
surface: with $m$ the multiplicity and $D$ as above, either $D$ is a
smooth elliptic curve (type $_ mI_0$) or $D$ is a reduced cycle of
$n$ rational curves (type $_ mI_n$).

\begin{defn}
Let $\pi:S\to B$ be an elliptic fibration. We say that the fibre
$F_p=\pi^{-1}(p)$, $p\in B$, is singular if it has positive Euler
number. In particular, a multiple fibre with smooth reduction is not
singular in this sense.
\end{defn}

Except the type $_ mI_0$, all other types  have positive Euler
numbers \cite{wall}, and are singular by our definition.

\begin{lemma}
A minimal elliptic surface $S$ has nonnegative Euler number $\chi$,
and  the case $\chi=0$ occurs if and only if $S$ has no singular
fibres.
\end{lemma}

\NI{\it Proof.} Let $\pi:S\to B$ be an elliptic fibration,
$\Xi\subset B$ be the set of critical values of $\pi$, and $F_b$ be
the fibre
over $b$, $b\in B$. \\
Set $F=\displaystyle\bigcup_{b\in\Xi} F_b$, which is a closed set in
$S$, and we have the exact sequence
$$
...\to H^i_c( S\setminus F, \R)\to H^i(S, \R)\to H^i(F, \R)\to
H^{i+1}_c(S \setminus F, \R)\to...
$$
where the subscript $c$ means the cohomology with compact support.
Therefore, we deduce the relation among the Euler numbers
$\chi(S)=\chi_c(S\setminus F)+\chi(F)$.\\
Now $\pi: S\setminus F\to B\setminus \Sigma$ is a topological
elliptic fibre bundle with fibre $E$, and hence we have
$\chi_c(S\setminus F)=\chi_c(B\setminus \Sigma)\chi(E)=0$. As a
result,
$$\chi(S)=\chi(F)=\sum_{b\in\Sigma}F_b,$$
and the lemma follows. \QED

\smallskip

\NI To each elliptic fibration, we can associate two fundamental
invariants \cite{bpv, koII} :
\begin{defn}

Let $\pi:S\to B$ be an elliptic fibration. We define the following
invariants associated to $S$:
\begin{enumerate}
\item The $j$-invariant (functional invariant) $j_S$ of $S$. Let
$U$ be the open subset of $B$ consisting of regular values of $\pi$.
Let $f_S$ be the holomorphic map which associate to each point $b\in
B$ the isomorphism class of the elliptic curve $\pi^{-1}(b)\in
\C^+/\p SL_2(\Z)$, where $\C^+$ is the upper half plane, and j be
the biholomorphic function $j: \C^+ /\p SL_2(\Z)\to \C$ induced by
the elliptic modular function $\tilde j:\C^+\to \C$. Let $j_S:=
j\circ f_S/1728: U\to \C$. By the stable reduction theorem
\cite{bpv}, $j_S$ has an extension to a holomorphic function from
$B$ to $\p_1$.

\item The homological invariant (global monodromy) of $S$, which is a sheaf
$G_S$ on the base $B$. Let $U$ be the open subset of $B$ consisting
of regular values of $\pi$ and $i:U\to B$ be the inclusion map. Then
we define $G_S:=i_*(R^1\pi_*\Z|U)$. The sheaf $R^1\pi_*\Z|U$ is
locally constant, and it can be constructed from a representation of
the fundamental group $L_S: \pi_1(U, b)\to Aut(H^1(\pi^{-1}(b),
\Z))\cong SL_2(\Z)$ .
\end{enumerate}
\end{defn}
These two invariants are not unrelated. There is a natural
compatibility between them. From the definition, the functional
invariant $j_S$ takes the value of $\infty$ only at singular fibres.
Let $\pi: S\to B$ be an elliptic fibration. Assume the functional
invariant $j_S$ is not identically 0 or 1.  Let $U:=B \setminus
j^{-1}(\{0,1, \infty\})$. Composing with the canonical projection
$SL_2(\Z)\to \p SL_2(\Z)$, the equivalence class of the
representation $G_S: \pi_1(U)\to  SL_2(\Z)$ induces a representation
$\tilde{G}_S:\pi_1(U)\to \p SL_2(\Z)$. On the other hand, the
elliptic modular function $\tilde j:\C^+\to \C$ is unbranched
outside the preimage of $0$ and $1$. The covering $\tilde
j:\C^+\setminus j^{-1}(\{0,1\})\to \C\setminus \{0,1\}$ therefore
induces an equivalence class of a representation $\tilde
j_*:\pi_1(\C\setminus \{0,1\})\to \p SL_2(\Z)$.  Then the
meromorphic function $j_S$ gives a map $j_{S*}: \pi_1(U)\to \p
SL_2(\Z)$. Making the identification of the fibre $F_b$ with the
quotient $\C/\Z\oplus\Z\tau$ before and after going around the loop,
we see that $j_{S*}$ can be thought of as the monodromy of the
period $\Z\oplus\Z\tau$. But it is slightly cruder than the
monodromy since it does not distinguish $\pm \tau$. The discussion
gives the following commutative diagram
$$
\begin{CD}
\pi_1(U) @>{\rm G_S}>> SL_2(\Z) \\
@V{j_{S*}}VV  @VVV\\
\p SL_2(\Z)@= \p SL_2(\Z).
\end{CD}
$$
\smallskip
\begin{defn}
Let $L:\pi_1(U)\to SL_2(\Z)$ be an equivalence class of a
representation, and $h:U\to\C$ be a meromorphic function with
$h(u)\neq 0,1,\infty$ for $u\in U$. We say $L$ belongs to $h$ if it
induces the representation defined by $h$.
\end{defn}


Kodaira has  the following result about the classification of
elliptic fibrations without multiple fibres.
\begin{thm}(Kodaira \cite{koII})\label{cl}
Let $B$ be a Riemann surface, $U= B \setminus \{p_1,...p_k\}$, and
$h$ be a meromorphic function on $B$ with $h(u)\neq 0,1,\infty$ for
$u\in U$.

\begin{enumerate}
\item If $k\geq 1$, then there are exactly $2^{g(B)+k-1}$
inequivalent homological invariants $L$ that belong to $h$, where
$g(B)$ is the genus of the curve $B$.
\item Given $h$ and a homological
invariant $L$ belonging to $h$, there exists a unique minimal
elliptic fibration $f:S\to B$ with these invariants, admitting a
section.
\item
Let $\mathscr{F}(h, L)$ denote the set of all elliptic fibrations,
without multiple fibres, with given invariants $h$ and $L$. Given
$S'$ an element of $\mathscr{F}(h, L)$, there exist a covering
$B=\cup V_i$ with $f _i: S_i:= f^{-1} (V_i)\to V_i$ being the
restriction of $f$, and a cocycle $\xi_{ij}\in H^1(\mathscr{T})$,
where $\mathscr{T}$ is the sheaf of local holomorphic sections of
the unique  elliptic surface in (2), such that $S'$ is obtained by
gluing the collections of $S_i$ together using the $\xi_{ij}$. In
particular, the set $\mathscr{F}(h, L)$ is parameterized by the
abelian group $H^1(\mathscr{T})$.
\end{enumerate}
\end{thm}

\begin{defn}
The unique minimal elliptic surface $S^{\mathfrak{B}}$, which admits
the invariants $j_S$ and $G_S$ of $\pi:S\to B$ and a section, is
called the basic elliptic surface associated with $S$.
\end{defn}

From Theorem \ref{cl}, we see that an elliptic surface without
multiple fibres is locally isomorphic to its basic elliptic surface.

\smallskip
\subsection{Elliptic fibre bundles}\

A special case of elliptic surfaces is the total space of a
holomorphic elliptic fibre bundle $\pi:S\to B$ which is locally
trivial. Let $E$ be a smooth elliptic curve. A holomorphic elliptic
fibre bundle $S \stackrel{\pi}\rightarrow B$ is determined by the
associated class $\xi\in H^1(B,{\mathscr{A}_B})$ where ${\mathscr
A}_B$ is the sheaf of germs of local holomorphic maps from $B$ to
$Aut(E)$.
\begin{defn}
Let $E$ be a smooth elliptic curve. The $E$-fibre bundle
$S\stackrel{\pi}\rightarrow B$ is called a principal bundle if the
structure group can be reduced to $E$.
\end{defn}
Let $S\stackrel{\pi}\to B$ be a principal $E$-bundle with the
associated class $\xi\in H^1(B, \mathscr E_B)$, where $\mathscr E_B$
is the sheaf of germs of local holomorphic maps from $B$ to $E$. The
long exact sequence induced by the universal covering sequence is
written as
$$H^1(B,\Gamma)\rightarrow H^1(B,{\mathcal O}_B) \rightarrow H^1(B,
{\mathscr E}_B) \stackrel{c}\rightarrow H^2(B,\Gamma) \rightarrow  0
,$$ where $\Gamma$ is the lattice such that $\C/\Gamma=E$.
Topologically, principal $E$-bundles are classified by the
associated Chern class $c(\xi)\in H^2(B, \Gamma)$. The Chern class
$c(\xi)$ vanishes if and only if the principal $E$-bundle $S\to B$
is topologically trivial, hence the first Betti number
$b_1(S)=b_1(B)+2$ is even. In this case, by chasing the following
diagram in cohomology
\[
\begin{CD}
H^1(B,\C) @>>>H^1(B, E) @>>>H^2(B, \Gamma)@>>> H^2(B,\C)\\
@VVV @VVV  \Big\Vert  &\\
H^1(B,\mathcal{O}_B) @>>>H^1(B,{\mathcal E}_B ) @>c>>H^2(B, \Gamma),& &\\
\end{CD}
\]
we find that the $E$-bundle is defined by a locally constant cocycle
in $H^1(B, E)$. On the other hand, if $c(\xi)\neq 0$, then the
splitting $E=S^1\times S^1$ would imply the existence of an
$S^1$-bundle $X$ over $B$ such that $S=X\times S^1$. Therefore the
Gynsin sequence
\[
\begin{CD}
0@>>>H^1(B,\Z) @>>>H^1(X, \Z) @>>>H^0(B, \Z) \\
@>\delta >>H^2(B,\Z) @>>>H^2(X,\Z) @>>>H^1(B, \Z)@>>> 0,\\
\end{CD}
\]
where $\delta$ is the multiplication with the Chern class $c(\xi)$,
tells us that $b_1(S)=b_1(B)+1$ is odd.

\begin{prop}\label{principal}
Let $E$ be a smooth elliptic curve, $\pi: S\to B$ be a principal
$E$-bundle. Then $S$ admits cscK metrics if and only if the first
Betti number $b_1$ of $S$ is even.
\end{prop}
\NI{\it Proof.} The only if part is an immediate consequence of the
Hodge decomposition for a compact K\"ahler manifold. To prove the if
part, the above discussion tells us that  the surface $S$ is defined
over some open covering $\{U_i\}$ of $B$ by patching the pieces
$U_i\times E$ together according to some locally constant cocycle in
$H^1(B, E)$. Endow the flat metric on $E$. Fix a cscK metric on the
curve $B$ (Depending on the genus of $B$, it is either the
Fubini-Study, flat, or hyperbolic metric.), and restrict it to each
open set $U_i$. Then we can see that the product metrics on all
pieces $U_i\times E$ are respected by the patching procedure since
the cocycle is locally constant. Therefore we obtain a globally
defined cscK metric on the surface $S$. \QED

In fact, we have $H^1(B, E)=Hom(H_1(B, \Z), E)=Hom(\pi_1(B), E)$
from the universal coefficient theorem. Let $\xi\in H^1(B, E)$ be
the locally constant cocycle which determines the principal
$E$-bundle $S\to B$. Regard $\xi$ as a homomorphism from $\pi_1(B)$
to $E$, which acts on $\mathbb{H}\times E$. Then the principal
$E$-bundle $S$ is obtained from the quotient of  $\mathbb{H}\times
E$ by the action of $\xi$. In other words, the universal cover of
every principal $E$-bundle with even first Betti number is
$\mathbb{H}\times \C$ and each deck transformation can be expressed
by $(z,w)\mapsto(\alpha(z),w+t)$ where $\alpha(z)\in \p SL_2(Z)$,
and $t$ is a constant.

\begin{cor}\label{trivial}
Let $\pi:S\to B$ be an elliptic fibre bundle with even first Betti
number. If the functional invariant $j_S$ is constant and the
homological invariant $G_S$ is trivial, then there exists an
elliptic curve $E$  such that $\pi:S\to B$ is a principal
$E$-bundle. In particular, $S$ admits a cscK metric.
\end{cor}
\NI{\it Proof.}  Since $S$ has no singular fibre and $j_S$ is
constant, all fibres are isomorphic to some elliptic curve $E$.
Since $j_S$ is constant and $G_S$ is trivial, the basic elliptic
surface of $S$ is the product $B\times E$. Using Theorem \ref{cl},
the surface $S$ is obtained from $B\times E$ by twisting according
to some cocycle $\xi\in H^1(B, \mathscr{E}_B)$. Therefore $S$ is a
principal $E$-bundle with even Betti number. The corollary follows
then from Proposition \ref{principal}. \QED

\begin{thm}\label{nomultiple}
Let $\pi:S\to B$ be an elliptic fibration   whose each fibre is
smooth and whose first Betti number $b_1$ is even. Assume that the
genus of $B$ is at least 2. Then there exists a principal $E$-bundle
$S'\to E$ which is a fibre-preserving \'{e}tale covering of $S\to
B$.
\end{thm}

\NI {\it Proof.} Since $S$ has no singular fibres, the $j$-invariant
$j_S$ can not take the value of $\infty$, therefore it is constant,
and every fibre $E$ is isomorphic.  Since $G_S$ belongs to $J_S$,
the homological invariant is equivalent to the monodromy
representation $G_S: \pi_1(B)\to \Z_m$. (Here the value of $m$
depends on the value of $j_S$: if $j_S$ is 0, then $m$ =2, 3, or 6;
if $j_S$ is 1, then $m$= 4, or 6; if $j_S$ is 1, then $m=2$ or 4;
otherwise, $m$ = 2.)  Therefore  there exists a connected unramified
cover $B'\to B$ of degree $m$ such that the pullback $S'= B'\times_B
S$ is an elliptic surface over $B'$ with fibre $E$, constant
$j_{S'}$, trivial $G_{S'}$, and even Betti number.  From Corollary
\ref{trivial}, $S'$ is obtained by twisting $B'\times E$ according
to some constant cocycle $\xi'\in H^1(B', E)$. Therefore $S'\to B'$
is a principal $E$-bundle and by construction,  it is an \'{e}tale
$m$-cover of $S$ and it is locally isomorphic to every fibre of
$S\to B$. The theorem follows. \QED

\smallskip
\subsection{Elliptic surfaces whose each fibre has a  smooth
reduction}\

One crucial tool in studying the elliptic surfaces is Kodaira's
formula (see \cite{bpv} p.161) of the canonical divisor of an
elliptic surface $\pi: S\to B$:
$$K_S= \pi^*(K_B+D)+\Sigma (m_i-1)F_i$$
where the $F_i$'s are the multiple fibres of multiplicity $m_i$ and
$D$ is some divisor of $B$ with $\mbox{deg}\, D = \chi(\mathcal
{O}_S)$. The formula implies that $K_S^2=0$. Usually we assume $B$
is an orbifold with orbifold points $P_i$ of order $m_i$
corresponding to a fibre $F_i$ of multiplicity $m_i$. Let
$\tau_S=\chi (\mathcal {O}_S)-\chi^ {orb}(B)$, where $\chi^{orb}(B)$
is the orbifold Euler number defined by
$$
\chi^{orb}(B):= \chi^{top}(B) - \sum_{j=1}^k (1- \frac 1 {m_j}).
$$
Using that the plurigenera $P_m(S)=h^0(S, K^{\otimes m})=m\tau_S$
when $m$ is divisible by $m_i$ for all $i$ and some extra thoughts,
Wall \cite{wall} shows that the sign of $\tau_S$ determines the
Kodaira dimension of an elliptic surface.

\begin{lemma}(Wall \cite{wall})\label{sign}
Let $\pi:S\to B$ be an elliptic fibration and $\tau_S$ be defined as
above. Then the Kodaira dimension $\kappa$ of S is $-\infty,0,$ or 1
corresponding to $\tau_S<0$, $\tau_S=0$, or $\tau_S>0$,
respectively.
\end{lemma}
\NI {\it Proof.} See \cite{wall} Lemma 7.1. \QED

In the following, we are interested in elliptic surfaces of Kodaira
dimension $\kappa=1$, which are sometimes called properly elliptic
surfaces. We will exploit some properties of properly elliptic
surfaces without singular fibres, and use them later to construct
elliptic surfaces with cscK metrics. Let us start with a feature of
the base of an elliptic fibration.

\begin{defn} An orbifold Riemann surface is called good if its
orbifold universal cover admits a cscK metric.
\end{defn}

\begin{lemma} \label{goodorbifold} Let $S$ be a minimal elliptic surface
with $\kappa(S)=1$ and without singular fibres. Then the base $B$ is
a good orbifold.
\end{lemma}
\NI {\it Proof.} Using Noether's formula, $\chi(\mathcal O_S)=1/12
(K_S^2+\chi)=1/12 \chi=0$. By Lemma \ref{sign}, we have
$\chi^{orb}(B)=-\tau_S<0$. The lemma follows from Troyanov's
argument \cite{tr} that an orbifold Riemann surface $B$ is always
good if $\chi^{orb}(B) \leq 0.$ \QED

Smooth elliptic fibration without multiple fibres are completely
classified by Kodaira \cite{koII}. If there are multiple fibres, we
need some more work to reduce it to a smooth fibration:

\begin{lemma}(local version)\label{local}
Let $\pi: S\to \bigtriangleup$ be an elliptic fibration over the
unit disk $\bigtriangleup$, and the fibre $F_0=\pi^{-1}(0)$ is a
multiple fibre with multiplicity $n$. Let $\delta_n:
\bigtriangleup\to\bigtriangleup$ be the map $z\to z^n$,
$\bar{S}=\bigtriangleup\times_{\bigtriangleup} S$ be the fibre
product with respect to $\delta_n$, and $S'$ be the normalisation of
$\bar S$. Then $S'$ is a nonsingular surface with no multiple fibre
and $S'\to S$ is an unramified cover.
\end{lemma}
\NI{\it Proof.} Let $x$, $y$ be the local coordinate on $S$ such
that $\pi(x,y)=x^n$. Then the fibre product $\bar S$ can be
expressed as
$$ \begin{array}{lll}
U\times_{\triangle}\triangle &=&\{(x,y,z)|x^n=z^n\}\\
&=& \bigcup\limits_{\zeta\in \mu_n }U_{\zeta},\\
\end{array}
$$
where $U_{\zeta}=\{(x,y,\zeta x)|(x,y)\in U\}\cong U$, and $\mu_n$
is the group of the n-th root of unity. The normalisation $S'$
locally is the disjoint union of the $n$ components $U_{\zeta}$, and
therefore is smooth. Identifying $U_{\zeta}$ with $U$, we can define
$\pi':S'\to\bigtriangleup$ by $(x,y) \mapsto \zeta x$ on
$U_{\zeta}$, which has no multiple fibre.  Moreover, the group
$\mu_n$ acts on $S'$ freely by interchanging these components
$U_{\zeta}$ such that $S'/u_n=S$.\QED

\begin{defn}
The fibration $\pi':S'\to \bigtriangleup$ constructed in the above
proof is called the $n$-th root fibration of $\pi$.
\end{defn}

\begin{lemma}\label{multiplesmooth}
(\cite{b} Lemma 6.7) Let $\pi: S\to B$ be an elliptic surface whose
base $B$ is a good orbifold and whose fibres are either smooth or
multiples of smooth elliptic curves. Then there exists a branched
Galois cover $p_1:B'\to B$ with Galois group $G$, a surface $S'$ and
a commutative diagram
$$\begin{CD}
  S' @>>> S   \\
  @VV\pi' V         @VV{\pi}V \\
  B' @>p_1>> B
\end{CD}$$
such that the action of $G$ on $B'$ lifts to $S'$, $p_1$ induces an
isomorphism $S'/G\to S$, and every fibre of $\pi'$ is smooth.
\end{lemma}
\NI{\it Proof.} Since $B$ is a good 2-orbifold Riemann surface,
there exists a smooth Riemann surface $B'$, and a branched covering
map  $p_1:B'\to B$ with the group of deck transformation $G$. We can
assume $p_1$ is a Galois cover. ( If it is not, then $\pi_1(B')$ is
not a normal subgroup of $\pi_1(B)$. We can then take the
intersection $H$ of all conjugate groups of $\pi_1(B')$ in
$\pi_1(B)$, which is a subgroup of $\pi_1(B')$ of finite index and
is normal in $\pi_1(B)$. Therefore there exists a further finite
cover $B''\to B'\to B$ such that the group of deck transformations
of $B''$ over $B$ is $H$, which acts transitively on the preimage of
any point on $B$.)  Now let $p_1:B'\to B$ be a Golois cover with the
Golois group $G$. We have $B'/G=B$. Consider the pull-back
$\bar{S}:=B'\times_{B} S$ which has singularities at the pull-back
of the multiple fibres. Let $\Sigma\subseteq B$ be the set of
multiple points. Outside the fibres over $\Sigma$, there is an
\'etale cover
$$
\bar{S}_1:=\bar{S}\setminus \bigcup\limits_{b\in
\Sigma}\pi'^{-1}(p_1^{-1} (b)) \longrightarrow
S_1:=S\setminus\bigcup\limits_{b\in \Sigma} \pi^{-1}(b),
$$
and the action of $G$ lifted to  $\bar{S}_1$ acts freely such that
$\bar{S_1}/G= S_1$. Let $S'$ be the normalisation of $\bar{S}$. From
Lemma \ref{local}, $S'$ is smooth,  and the action of $G$ lifted to
$\bar S'$  acts on $S'$ freely
 such that $S'/G=S$.\QED

\begin{thm}\label{hc}
Let $S$ be  a minimal properly elliptic surface of K\"{a}hler type
whose fibres are either smooth or multiples of smooth
elliptic curves. Then the universal cover
is biholomorphic to $\mathbb{H} \times \C$, and $S$ inherits a cscK metric from $\mathbb{H} \times \C$.
\end{thm}

\NI{\it Proof.} (cf.\cite{wall}) Let $\pi: S\to B$ be the elliptic
fibration. By Lemma \ref{multiplesmooth}, there exists a smooth
finite cover $B'$, and a branched Golois covering map $p_1:B'\to B$
with Galois group $G$. Lift the action of $G$ to the normalization
$S'$ of the pullback $B'\times_{B} S$, which is a smooth elliptic
fibration over $B'$. Then we can see that $G$ acts freely on $S'$
and the quotient $S'/G=S$. From Theorem \ref{nomultiple}, there is
an  \'{e}tale cover $B''$ of $B'$   such that the pull-back surface
$S''$ is a principal $E$-bundle with even first Betti number. That
is,  $S''$ is obtained by twisting the trivial bundle $B''\times E$
by a locally constant cocycle in $H^1(B'', E)$. We can then see
$\mathbb{H}\times \C$ is the universal cover of $S'', S'$, and $S$.
And the action of every deck transformation $\gamma$ of
$\mathbb{H}\times \C$ over $S''$ can be written explicitly as
$$\gamma(z,w)=(\alpha_{\gamma}(z), w+t_{\gamma}),$$
where $\alpha_{\gamma}\in \p SL_2(\C)$, $t_{\gamma}\in \C.$ So far,
we have the following diagram
\[
\begin{CD}
\mathbb{H}\times \C @>>> S'' @>>> S\\
@VVV @VVV @VVV\\
\mathbb{H} @>>> B''@>>> B
\end{CD}
\]

Without loss of generality, we can assume $S''$ is a normal cover of
$S$. Let $H$ be the group of deck transformations of
$\mathbb{H}\times \C$ over $S$, and $\Gamma$ be the group of deck
transformations of $\mathbb{H}\times \C$ over $S''$. The action of
$\Gamma$ on $\mathbb{H}\times\C$ descends to a discrete group
$\bar\Gamma$ acting on the base $\mathbb{H}$ and the quotient
$\mathbb{H}/\bar\Sigma=B''$ is the base of $S''$. By theorem
\ref{principal}, we have that the group $\Gamma$ is a subgroup of
$\p SL_2\times \C$ and $S''$ inherits the product metric from
$\mathbb{H}\times \C$.  To show $S$ inherits the product metric from
$\mathbb{H}\times \C$, it suffices to show that every deck
transformation  $h\in H$  preserves the product metric. Given $h\in
H$, since all fibres of $S''$ are isomorphic, and the group of deck
transformations of $S''$ over $S$ preserves the elliptic structure,
we can write the action of $h\in H$ on $\mathbb{H}\times \C$ as
$$h(z,w)=(\alpha_h(z), \epsilon_h w+a_h(z)),$$
where $\epsilon_h$ is the $k$-th root of unity for $k = 2,4,\mbox{or
} 6$. For a covering transformation $h$, $\epsilon_h$ is constant. A
direct computation shows that
\smallskip
\begin{align}
h^{-1}(z,w)&=(\alpha^{-1}_h(z), \epsilon^{-1}_h w-\epsilon^{-1}_ha_h(\alpha^{-1}_h(z))),\nonumber\\
h\gamma h^{-1}(z,w)&=(\alpha_h\alpha_{\gamma}\alpha_{h^{-1}}(z),
w-a_h(\alpha_{h^{-1}}(z))+\epsilon_h
t_{\gamma}+a_h(\alpha_{\gamma}\alpha_{h^{-1}}(z)))\nonumber.
\nonumber
\end{align}
It follows that $a_h(\alpha_{\gamma}(z))-a_h(z)=t_{h\gamma
h^{-1}}-\epsilon_ht_{\gamma}$ is a constant. For a fixed $h$, the
cocycle $\{t_{h\gamma h^{-1}}-\epsilon_ht_{\gamma}\}\in
Hom(\pi_1(B''), \C)\cong H^1(B'',\C)$ defines a principal
$\C$-bundle $X$ over $B''$, for which $a_h(z)$ defines a section.
Therefore the principle $\C$-bundle is trivial, and $t_{h\gamma
h^{-1}}-\epsilon_ht_{\gamma}=0$. It follows that the holomorphic
function $a_h(z)$ is constant on the orbits of $\bar\Gamma$, and can
be regarded as a holomorphic function of the compact surface $B''$.
Thus $a_h(z)$ is constant. This concludes the proof.\QED

\subsection{Logarithmic transform}\

Logarithmic transform is an important operation introduced by
Kodaira \cite{koII}. It enables us to replace a smooth fibre in an
elliptic fibration by a multiple fibre. (For more details see
\cite{gh}, Chap. 4.5). Let $\pi:S\to B$ be an elliptic fibration
over a smooth curve $B$. Let $U \subset B$ be an open neighborhood
of $p \in B$, $E_0$ be the smooth fiber over $p$. Choose some local
coordinate $z$ with $z(0)=p$. Denote $\Sigma:=\pi^{-1}(U)$. For
every $m\in \mathbb N$, consider the diagram

$$\begin{CD}
  \tilde{\Sigma} @>>> \Sigma   \\
  @VV{\tilde{\pi}}V         @VV{\pi}V \\
  \tilde {U} @>\phi>> U
\end{CD}$$
where $\phi$ is a cyclic cover of degree $m$ given by $
\tilde{z}\mapsto \tilde{z}^m$, and $\tilde{\Sigma} =\phi^*
\Sigma=\tilde U\times_{U}\Sigma$. Let $\Lambda(z)$ be a family of
lattices such that $\Sigma=U\times\mathbb{C}/\Lambda(z)$, then
$\tilde\Sigma=\tilde{U}\times\mathbb{C}/\Lambda(\tilde{z}^m)$. Let
$\beta(z)$ be a local $m$-torsion section of $\Sigma\to U$ and $G\in
Aut(\tilde{\Sigma})$ be the cyclic group  generated by
$$
(\tilde{z}, t)\mapsto (e^{\frac{2\pi\sqrt{-1}}{m}}\tilde{z},
t+\beta(\tilde{z}^m)) \,\,\mod \Lambda(\tilde{z}^m).$$ Denote
$\Sigma':=\tilde{\Sigma}/G$. We then have an isomorphism
$\alpha:\Sigma'\backslash(\phi^*E_0/G)\cong \Sigma \backslash E_0$,
where $E_0=\pi^{-1}(p)$. By setting $S'=(S\backslash
E_0)\bigcup_{\alpha} \Sigma'$, we get an elliptic fibration $S'\to
B$, which is isomorphic to $S$ away from $E_0$ and has a multiple
fibre with reduction isomorphic to $E_0/G$ of multiplicity $m$ over
$p$.

Since we simply remove $T^2\times D^2 (\approx \Sigma)$ and glue it
back in (as $\Sigma'$) with a different fibration, we can
reformulate the logarithmic transformation in the following way:
Consider an elliptic surface $\pi:S\to C$ and fix a generic fibre
$\pi^{-1}(t)=F$. We denote a closed tubular neighborhood of the
fibre $F$ in $S$ by $\bar{\Sigma}$ with the interior $\Sigma$;
$\bar{\Sigma}$ is diffeomorphic to $T^2\times D^2$. Deleting the
interior $\Sigma$ from $S$ and regluing $T^2\times D^2$ via a smooth
map $\phi:T^2\times S^1 \to$ Boundary$(S \setminus \Sigma)$ with
multiplicity $m$, we get a new manifold $S'$. The diffeomorphism
type of the resulting manifolds depends on the multiplicity of $m$
of the map $\phi$.

A remark here is that the process of logarithmic transform is quite
violent. A non-algebraic surface may be obtained from an algebraic
one by a logarithmic transform, and vice versa.

\section{A construction of elliptic surfaces of cscK metrics}
In this section, we construct elliptic surfaces with cscK metrics.

\subsection{Arezzo-Pacard Theorem and some results}\

The main tool we will use  is  the following remarkable theorems
\ref{ap1}, \ref{ap2} by Arezzo-Pacard \cite{ap}. Kronheimer
\cite{kr} has shown that there exist asymptotically locally
Euclidean resolutions of singularities of the type $\C^2/\Gamma$,
where $\Gamma$ is a finite subgroup of $SU(2)$. Arezzo and Pacard
\cite{ap} glue this resolution to an orbifold with isolated
singularities of local model $\C^2/\Gamma$, which has no nontrivial
holomorphic vector fields, study  the partial differential equations
arising from the perturbation of the K\"ahler forms suggested in
\cite{ls}, and they succeed constructing cscK metrics on the
resulting desingularisation. This is a major breakthrough in the
construction of manifolds with cscK metrics, and we state their
results here.

\begin{thm}(Arezzo, Pacard \cite{ap}) \label{ap1} Let $M$ be a compact
cscK manifold. Assume there are no nontrivial  holomorphic vector
fields that vanish somewhere on $M$. Then, given finitely many
points $p_1, p_2,..., p_n$ on $M$, the blow-up of $M$ at $p_1,
p_2,..., p_n$ carries a cscK metric.
\end{thm}
\begin{thm}(Arezzo, Pacard \cite{ap}) \label{ap2} Let $M$ be a compact
$2$-dimensional cscK orbifold with
isolated singularities. Assume there are no
nontrivial holomorphic vector fields on $M$.  Let $p_1, p_2,...,
p_n$ be  finitely many points on $M$, each of which has a
neighborhood biholomorphic to $\C^2/\Gamma_j$, where $\Gamma_j$ is a
finite subgroup of $SU(2)$. Then the minimal resolution of $M$ at
$p_1, p_2,..., p_n$ carries cscK metrics.
\end{thm}

The first result we obtain by applying Arezzo-Pacard Theorem \ref{ap1},
\ref{ap2} is that
an elliptic surface obtained by blowing up a minimal properly
elliptic surface of K\"ahler type, which has no singular fibres, at
finitely many points admits cscK metrics.

\begin{prop}\label{parallel} Let $S$ be a
minimal properly elliptic surface of K\"{a}hler type, which has no
singular fibres. Then every nontrivial holomorphic vector fields on
$S$ has no zeros.
\end{prop}

\NI{\it Proof.} By Lemma \ref{goodorbifold}, the assumption implies
the base $B$ is a good orbifold. Let $\mathfrak h(S)$ be the algebra
of holomorphic vector fields. Let $\zeta \in\mathfrak h(S)$ be
nontrivial. Since the base $B$ is a 2-orbifold with
$\chi^{orb}(B)<0$, $\zeta$ is vertical. By Theorem \ref{hc},
$\mathbb H\times \C$ is the universal cover of $S$. Thus $\xi$ can
be lifted to a holomorphic vector field $\tilde{\zeta}=
f(z)\frac{\partial}{\partial w}$ on $\mathbb H \times \C$, where $f$
is a holomorphic function on $\mathbb H$, $z$ is the coordinate of
$\mathbb H$, and $w$ is the coordinate on $\C$. Since $f$ is
invariant under the action of $\pi_1(B)$, it defines a holomorphic
function on $B$, and is constant. Therefore we conclude that $\zeta$
is induced by $c\frac{\partial}{\partial w}$, which is everywhere
nonzero. \QED

\begin{cor}
Let $S$ be a minimal properly elliptic surface of K\"{a}hler type,
which has no singular fibres. Then every surface obtained by blowing
up $S$ at finitely many points admits cscK metrics.
\end{cor}
\NI{\it Proof.} The corollary follows from Arezzo-Pacard Theorem
\ref{ap1} and Proposition \ref{parallel}. \QED

\subsection{The construction of elliptic surfaces with cscK
metrics}\

Recall that for a minimal elliptic surface $S$,
 using Noether's formula, we have $\chi(\mathcal O_S)=1/12
(K_S^2+\chi)=1/12 \chi$. Therefore, the Euler number $\chi(S)$ is
always a multiple of 12.

\begin{lemma}
Let $\chi$ be a positive multiple of 12, $g$ be a nonnegative
integer. Then there exists a smooth compact Riemann surface
$(\Sigma_{g'}, J)$ of genus $g':=\frac{1}{12}\chi+2g-1$, which
admits an isometric holomorphic involution $\psi$  with
$\frac{1}{6}\chi$ fixed points.
\end{lemma}
\NI{\it Proof.} When $g'$ is at least 2, such a surface
$(\Sigma_{g'}, J, \psi)$ exists due to Thurston's pants
decomposition \cite{p} of a Riemann surface $\Sigma_{g'}$, which
states that a hyperbolic Riemann surface $R$ of genus $l$ always
contains $3l-3$ simple closed geodesics such that cutting $R$ along
these geodesics decomposes $R$ into $2g-2$ pairs of pants, and the
fact that for any pair of pants, the length of boundaries  provide
the coordinates for Teichm\"uller space, and can be chosen
arbitrarily. When $g'=0,1$, this involution $\psi$ amounts to a
$180^{\circ}$ rotation of $S^2$ about an axis, or the Weierstrass
involution of a torus.\QED

This lemma tells us that given two numbers $\chi=12d>0$, and $g \geq
0$, there exists a Riemann surface of genus
$g'=\frac{1}{12}\chi+g-1$, and it admits a holomorphic surjective
holomorphic map $\Sigma_{g'}\to B$ of degree 2, where $B$ is a
Riemann surface of genus $g$.

\smallskip

Our goal is to construct  an elliptic surface with Euler number
$\chi$, base curve $B$ of genus $g$ and $k$ multiple fibres  of
multiplicities $m_1, m_2,...,m_k$, respectively, where the numbers
$\chi$, $g$, $k$ and $m_1, m_2,..., m_k$ are given.

\begin{figure}[htbp]
  \begin{center}
    \mbox{
      \subfigure[$g \geq 1, \chi\geq 24$]{\scalebox{0.75}{\input{1.pstex_t}}} \quad
      \subfigure[$g \geq 1, \chi= 12$]{\scalebox{0.85}{\input{2.pstex_t}}} \quad
      }
  \end{center}
\end{figure}
\begin{figure}[htbp]
  \begin{center}
    \mbox{
      \subfigure[$g=0, \chi\geq 24$]{\scalebox{0.75}{\input{3.pstex_t}}}\quad
      \subfigure[$g=0, \chi=12$]{\scalebox{0.85}{\input{4.pstex_t}}} \quad
      }
  \end{center}
\end{figure}


Now let $E=\C/\Gamma$ be a fixed elliptic curve, and $\phi :
E\rightarrow E$ be the Weierstrass involution, which is an isometry
with respect to the flat metric on $E$ with 4 fixed points. Let
$\hat{S}$ be the product $\Sigma_{g'} \times E$. Choose arbitrarily
$k$ points $P_1, P_2,..,P_k$ on $\Sigma_{g'}$ outside the fixed
points of $\Psi$. Let $Q_i:=\psi(P_i)\in \Sigma_{g'}$. Apply the
logarithmic transformation of order $m_i=n_i$ at theses $2k$ points
$P_i$, $Q_i$, where $i=1,...,k$. Let $S$ be the resulting surface.
The involution $\psi \times \phi$ on the product $\Sigma_{g'} \times
E$ can be extended to an involution  $\iota$ on $S$, which has
$\frac{2}{3}\chi$ fixed points.

From the construction, $S\to\Sigma_{g'}$ is an elliptic fibration
with $2k$ multiple fibres.  Now take the quotient of $S$ by the
action of $\iota$. The resulting surface $S'= S /\Z_2$ is singular
and has $\frac{2}{3}\chi$ ordinary double points. Let $S''$ be a
minimal resolution, which is obtained by replacing each singular
point by a $(-2)$-curve. we can see that $S''\to B$ is an elliptic
fibration with  $k$ multiple fibres of order $m_1, ..., m_k$
respectively.

\begin{prop}
The elliptic surface $S''$ constructed above has Euler
characteristic number $\chi$.
\end{prop}

\NI{\it Proof.} Recall from the construction that the
$\frac{2}{3}\chi$ double points come from the fixed points of
$\iota$, and they are locally modeled by $\C^2/\Z_2$. To resolve the
singularities, we take the blow-up $\tilde S$ of $S$ at the
$\frac{2}{3}\chi$ double points, and extend the map $\iota$ to an
automorphism $\tilde{\iota}$ of $\tilde S$, then $S''= \tilde S /
\Z_2$ is the nonsingular complex surface obtained by replacing each
double point with a ($-2$)-curve.
Let $p_1:\tilde{S} \rightarrow S''$, and  $p_2:\tilde{S} \rightarrow
S$ be the quotient map and the blow-down, respectively.

\[
\begin{CD}
S @<p_2<<\tilde{S} @>p_1>>S''\\
@VVV @VVV  @VVV  \\
\Sigma_{g'}@= \Sigma_{g'}@>>>B\\
\end{CD}
\]

\NI We claim that the irregularity $q:=h^0(S'', \Omega^1_{S''})= g$
: If $\eta$ is a nonzero holomorphic 1-form on $S''$, then
$p_1^*{\eta}$ is  a holomorphic 1-form on $\tilde {S}$ invariant
under $\tilde{\iota}$. Since $H^0(S, \Omega^1_S) \rightarrow
H^0(\tilde{S}, \Omega^1_{\tilde{S}})$ is an isomorphism by Hartog's
extension theorem, there exists a 1-form $\xi$ on $S$ invariant
under $\iota$, such that $ p_2^{*}(\xi)=p_1^{*}(\eta)$. Since $S$ is
an elliptic surface without singular fibres, its universal cover is
$\mathbb H\times \C$, and we know that every holomorphic 1-form on
$S$ comes from either the base $\Sigma_{g'}$ or the fibre. However,
since $\phi^*dz=-dz$ on the elliptic curve $E$, the only 1-forms on
$S$ invariant under $\iota$ are induced from holomorphic 1-forms of
$\Sigma_{g'}$ invariant under $\psi$. That is, $q=\mbox{dim}_{\C}\{
\omega \in H^0(\Sigma_{g'},\Omega^1) |\psi^*\omega=\omega\}$. Since
every holomorphic 1-form on $\Sigma_{g'}$ invariant under the action
of $\psi$ corresponds to a holomorphic 1-form on
$\Sigma_{g'}/\psi=B$, we conclude that $q=g$. Using similar argument
as above, we can also show $p_g=1/12\chi+1 +g$ : the only
holomorphic 2-forms on $S$ invariant under $\iota$ are induced from
the exterior product of holomorphic 1-forms of $E$ and those of
$\Sigma_{g'}$. Because $\phi^*dz=-dz$, we can take "anti-symmetric"
1-forms on $\Sigma_{g'}$, i.e.  $\psi^*\omega=-\omega$. Since every
1-form on $S$ can be decomposed as the sum of a symmetric 1-form and
an anti-symmetric one, we have $p_g=h^0(\Sigma_{g'},
\Omega^1)-g=(\frac{1}{12}\chi+2g-1)-g = \frac{1}{12}\chi +g-1.$ From
Noether's formula, we can therefore compute the Euler number of
$S''$ as
$$
\begin{array}{ccl}
12\chi(\mathcal{O}_{S''})-K_{S''}^2&=&12(h^0(\mathcal{O}_{S''})-
h^1(\mathcal{O}_{S''})+h^2(\mathcal{O}_{S''}))-0\\
&=&12(1-g+(\frac{1}{12}\chi+g-1))\\
&=&\chi\,  ,
\end{array}
$$ where the first equality holds since $S''$ is minimal.
\QED

\begin{prop}\label{even}
Let $S''$ be an elliptic surface constructed above. If $S''$ has
even first Betti number $b_1$, then $S''$ admits cscK metrics.
\end{prop}

\NI{\it Proof.} Since $S''$ has even first Betti number, it admits a
K\"ahler form $\omega$. Note that the pullback form $p_1^*\omega$ is
positive semi-definite on $\tilde S$, which is degenerate on the
exceptional divisors.  Consider the Eguchi-Hanson metric on the
total space of $\mathcal{O}(-1)$. By gluing the K\"ahler potentials
carefully, we can get a K\"ahler form $\tilde\omega$ on $\tilde S$,
therefore the first Betti number $b_1(\tilde S)$ is even. Because
blowing up at points does not change the first Betti number,
$b_1(S)$ is even too. Therefore $S$ is a minimal properly elliptic
surface of K\"ahler type, which has no singular fibres. By Theorem
\ref{hc}, $S$ admits cscK metrics. Furthermore, the involution
$\iota:S\rightarrow S$ is an isometric automorphism. This gives
$S'=S/\iota$ an cscK orbifold metric. In view of Arezzo-Pacard
Theorem \ref{ap2}, to show that the minimal resolution $S''$ admits
a cscK metric, it suffices to show that there are no nontrivial
holomorphic vector fields on the orbifold $S'$: If $\zeta'$ is a
holomorphic vector field on $S'$, it can be lifted to a holomorphic
vector field $\zeta$ on $S$, which is invariant under the action of
$\iota$, and $\zeta$ has zeros at the fixed points of $\iota$. By
Proposition \ref{parallel}, $\mathfrak h(S)$ consists of only
parallel holomorphic vector fields. Therefore $\zeta$ is trivial,
and so is $\zeta'$. \QED

\begin{prop}\label{ellipticpositive}
Let $S''$ be an elliptic surface constructed as above. If $S''$ has
even first Betti number $b_1$, then $S''$ admits no nontrivial
holomorphic vector fields. In particular, every blow-up of $S''$ at
finitely many points admits a cscK metric.
\end{prop}

\NI{\it Proof.} Let $\zeta \in\mathfrak h(S'')$ be a holomorphic
vector field.  The vector field $\zeta$ restricted to the
$(-2)$-curves is tangent to these curves, therefore $\zeta$ can be
lifted to a holomorphic vector field $\tilde\zeta$ on $\tilde S$,
which is invariant under $\tilde \iota$. Let $p_2: \tilde S \to S$
be the natural blow-down map. By Hartog's extension theorem, the
push-forward $p_{2*}\tilde\zeta$ is well-defined holomorphic vector
field on $S$ which vanishes at all fixed points of $\iota$. As we
have shown in the proof of Theorem \ref{even}, $b_1(S)$ is even.
Therefore, $S$ is a K\"ahler type minimal elliptic surface with
$\kappa(S)=1$ and has no singular fibres. By Proposition
\ref{parallel}, $p_{2*}\tilde\zeta$ can only  be parallel, and therefore
it vanishes everywhere.  The last statement follows from
a direct application of Arezzo-Pacard Theorem \ref{ap1}. \QED

\section{compact complex surfaces and cscK metrics}
In this section, we will show that every complex surface, which is
not in the deformation class of $\p_2\#k\overline{\p}_2$, $k=1$ or
2, with $b_1$ even is deformation equivalent to a surface which
admits cscK metrics. We will proceed with the help of classification
theorem of compact complex surfaces: Let $S$ be a complex surface,
and $K_S$ be the canonical line bundle. We can define the
pluri-canonical map $\iota_{K_S^{\otimes k}}:S \to \p(H^0(S,
K_S^{\otimes k}))^*$, which is a rational map, not defined at the
base locus of the linear system $|K_S^{\otimes k}|$. The Kodaira
dimension $\kappa$ of $S$ is defined to be the maximal dimension of
the image $\iota_{K_S^{\otimes k}}(S)$ for $k\geq 1$. Recall the
definition of deformation equivalence first:

\begin{definition}
Two smooth compact complex manifolds $M, N$ are said to be
deformation equivalent or of the same deformation type if there
exist connected reduced complex spaces $\mathcal{X}$ and $\mathcal
{T}$, and a proper holomorphic submersion
$\Phi:\mathcal{X}\to\mathcal{T}$, together with points $t_1, t_2
\in\mathcal{T}$ such that for each $t\in \mathcal{T}$,
$M_t:=\Phi^{-1}(t)$ is a compact complex submanifold, and
$\Phi^{-1}(t_1)=M$, $\Phi^{-1}(t_2)=N$.
\end{definition}

Since we consider only the reduced spaces, an equivalent definition
would be to assume that $\mathcal{T}$ consists of finitely many
irreducible components, each of which is smooth (and which can be
taken to be a disk in $\C$). Let $S_1$ and $S_2$ be two deformation
equivalent surfaces, and let $\tilde{S_1}$ and $\tilde{S_2}$ be the
blow-ups of $S_1$ and $S_2$ at $r$ points. A straightforward
argument shows that  $\tilde{S_1}$ and $\tilde{S_2}$ are again
deformation equivalent.

\subsection{Kodaira dimension $\kappa=2$} \

A complex surface $S$ in this case is said to be of general type.
Since we have $c_1^2(K)>0$, every minimal surface of general type is
projective. For a minimal surface of general type,
$\iota_{K_S^{\otimes k}}$ is a globally defined map for $k\geq 5$,
and it is an embedding away from some smooth rational $-2$ curves.
The image of these curves are isolated singular points with local
structure group $\Gamma$, where $\Gamma$ is a finite subgroup of
$SU_2$ (see \cite{bpv}). One can get the pluricanonical model
$X=\iota_{K_S^{\otimes k}}(S)$ by collapsing these $(-2)$-curves. If
$M$ has no $(-2)$-curves, it has negative first Chern class
$c_1(M)$, and Aubin-Yau Theorem \cite{au,yau} asserts that every
manifold with negative first Chern class admits a K\"ahler-Einstein
metric. Otherwise, Kobayashi \cite {ko} has shown that the
pluricanonical model $X$ admits a K\"ahler-Einstein orbifold metric
of negative scalar curvature by extending Aubin's proof of the
Calabi conjecture. Along with the fact that a complex manifold of
general type has no nontrivial holomorphic vector fields, a direct
application of Arezzo-Pacard Theorem \ref{ap1} gives the following
result.
\begin{thm}(Arezzo-Pacard \cite{ap})\label{general}
Every compact complex surface of general type admits cscK metrics.
\end{thm}
\NI{\it proof.} See \cite{ap} Corollary 8.3. \QED

\subsection{Kodaira dimension $\kappa=1$} \

Every complex surface $S$ in this case is a properly elliptic
surface. Since an elliptic curve has Euler number $\chi =0$, the
Euler characteristic number of an elliptic surface $S$ is given by
the sum of those of singular fibres. In particular, $\chi(S)\geq 0$
and the equality holds if and only if there are no singular fibres.

Now we introduce the following classification of the deformation
types of elliptic surfaces with singular fibres.

\begin{thm}\label{deformationelliptic}
Two elliptic surfaces with positive Euler numbers are deformation
equivalent (through elliptic surfaces) if and only if they have the same Euler
numbers and their base orbifolds are diffeomorphic.
\end{thm}
\NI{\it Proof.} See \cite{fm} Chap.1 Theorem 7.6 .\QED

\begin{thm}
A properly elliptic surface of K\"ahler type is deformation
equivalent to a compact complex surface with cscK metrics.
\end{thm}

\NI{\it Proof.} Let $\check{S}$ be the minimal model of $S$. If
$\chi(\check{S}) =0$, then $\check{S}$ is a minimal properly
elliptic surface of K\"ahler type which has no singular fibres. This
is done in  Theorem \ref{parallel}. If $\chi(\check{S})>0$,
a direct
application of Theorem \ref{deformationelliptic} shows that
$\check{S}$ is deformation equivalent to one of the elliptic
surfaces we constructed in Section 3. Therefore by Proposition
\ref{ellipticpositive}, $S$ admits cscK metrics. \QED

\subsection{Kodaira dimension {$\kappa=0$}}\

Minimal compact complex surfaces of K\"ahler type with Kodaira
dimension $\kappa=0$ consist of Enrique surfaces, K3 surfaces,
bielliptic surfaces, and Abelian surfaces. Although not all of them
are projective, they all admits K\"ahler metrics. Moreover, each of
them has a vanishing real first Chern class. A direct application of
the following celebrated theorem by  Yau \cite{yau} implies that
every compact K\"ahler manifold with vanishing real first Chern
class admits a Ricci-flat K\"ahler-Einstein metric.

\begin{thm}(Yau \cite{yau})
Let $M$ be a compact K\"ahler manifold with  K\"ahler form $\omega$,
and $c_1(M)$ be its real first Chern class. Then every closed real
2-form of (1,1)-type belonging to the class $2\pi c_1(M)$ is the
Ricci form of one and only one K\"ahler metric in the K\"ahler class
$[\omega]$.
\end{thm}

Since every holomorphic vector field on a Ricci-flat
K\"ahler manifold is parallel, there is no obstruction on the application of
Arezzo-Pacard Theorem \ref{ap2}, and we can reach the following result.

\begin{thm}\label{tori}
Let $S$ be a compact complex surface of K\"ahler type and of Kodaira
dimension $\kappa=0$. Then $S$ admits cscK metrics.
\end{thm}

\subsection{Kodaira dimension {$\kappa=-\infty$}}\

In this case, $S$ is called  a ruled surface. A minimal ruled
surface is either a geometrically ruled surface or the complex
projective plane $\p ^2$. If $S$ is a geometrically ruled surface,
then there exists a holomorphic rank two vector bundle $V$ over a
curve $C$ such that $S$ is isomorphic to the associated
$\p^1$-bundle $\p(V)$. Two vector bundles $V$, $V'$ over $C$ give
isomorphic ruled surfaces if and only if $V'=V\otimes L$ for some
holomorphic line bundle $L$ over $C$. It follows that $c_1(V)$ mod 2
is a holomorphic invariant of $S$. In addition, given rank two
vector bundles $V$ and $V'$ over the same curve $C$, if
$c_1(V)\equiv c_1(V')$ mod 2, then the resulting surfaces $\p(V)$
and $\p(V')$ are deformation equivalent \cite{fm}. Let $\pi:
S=\p(V)\rightarrow C$ be the ruling, $\sigma_0$ be the class of a
holomorphic section of $\pi$ (it always exists!), and $f$ be the
class of a fibre of $\pi$. We can see $\{\sigma_0, f\}$ is a basis
of $H^2(S,\Z)$ and $\sigma_0^2\equiv c_1(V)$ mod 2. Moreover, the
intersection pairing on $H^2(S,\Z)$ is even if $c_1(V)\equiv 0$ mod
2 and is odd if $c_1(V)\equiv 1$ mod 2. In particular, the
deformation equivalence class of a geometrically ruled surface
$\pi:S=\p (V)\to \C$ is determined by  the genus $g(C)$ of $C$ and
$c_1(V)$.

For every Riemann surface $C$ of genus $g\geq 2$, we shall show that
in each deformation class of geometrically ruled surfaces $\pi: S\to
C$, there exists one which admits cscK metrics. The justification
splits into two parts:

\subsubsection{$c_1(V)=0$ mod 2 }\

\label{representation} Let $C$ be a Riemann surface with the
hyperbolic metric, $\pi_1(C):=\langle a_1, b_1,...,a_g,
b_g:[a_1,b_1][a_2, b_2]...[a_g, b_g]=1 \rangle$ be the fundamental
group of $C$, and $g^{FS}$ be the (multiple of) Fubini-Study metric
with scalar curvature 1 on the complex projective line $\p_1$.
Define the representation $\rho:\pi_1(C)\rightarrow SU(2)/Z_2$ by

$$
\rho(a_1)=\rho(b_1):= \left[\left( \begin{array}{cc}
i &  0  \\
0 &  -i
\end{array} \right)\right],
$$

$$ \rho(a_2)= \rho(b_2):=\left[\left( \begin{array}{cc}
0 &  i  \\
i &  0
\end{array} \right)\right],
$$

\NI and $\rho(a_j)= \rho(b_j):=[id]$ for $j\neq 1,2$. Then
$\pi_1(C)$ acts isometrically on $\mathbb H \times  \p^1$, and the
quotient
$$S:=\mathbb{H}\times_{\rho} \p^1 =\mathbb{H}\times  \p^1/ \pi_1(C)$$
inherits a K\"{a}hler metric $g$ of constant scalar curvature
$s=-1+1=0$.

\begin{prop}\label{quasistable}
The surface $S$ defined above is deformation equivalent to the
trivial product $C\times \p^1$, and it has no nontrivial holomorphic
vector fields.
\end{prop}

\NI{\it Proof.} To prove the first statement, we construct a family
of surfaces defined by the following representations
$\rho_t:\pi_1(C)\rightarrow SU(2)/Z_2$:

$$
\begin{array}{l}
$$
\rho_t(a_1)=\rho_t(b_1)= \left[\left( \begin{array}{cc}
e^{it} &  0  \\
0 &  e^{-it}
\end{array} \right)\right],
$$\\

$$ \rho_t(a_2)= \rho_t(b_2)=\left[\left( \begin{array}{cc}
cos t &  e^{it} sin t \\
-e^{-it}sint &  cos t
\end{array} \right)\right],
$$
\end{array}
$$

\NI and $\rho_t(a_j)= \rho_t(b_j)=[id]$ for $j\neq 1,2$. Let
$S_t=\mathbb{H} \times_{\rho_t} \p^1$. This defines a family of
geometrically ruled surfaces in the same deformation class with
$S_0=C\times \p^1$ and $S_{\pi/2}=S$. To show the algebra $\mathfrak
h(S)$ of holomorphic vector fields consists of only the zero vector
field, we first notice that  the action of $\pi_1(C)$ by $\rho$
fixes no points of $\p^1$ by direct computation. Since a Riemann
surface with genus $g \geq 2$ admits no nontrivial holomorphic
vector fields, every holomorphic vector field $\xi$ on $S$ is
vertical. In particular, $\xi$ has zeros and is nonparallel. Since
$S$ carries a scalar flat K\"{a}hler metric, the algebra $\mathfrak
h_0$ of nonparallel holomorphic vector fields is the
complexification of the Lie algebra of nonparallel killing fields.
Lift $\xi$ to a vertical holomorphic vector field  $\tilde\xi$ on
$\mathbb H \times \C\p^1$, so $\tilde\xi \in \mathfrak {su}_2$,
which is invariant under the action of $\pi_1(C)$. This action of
$\pi_1(C)$ is given by composing the representation $\pi_1(C)\to
SU_2/Z_2$ with the adjoint action of $SU_2/Z_2$ on its Lie algebra.
Since the adjoint action of $SU_2/\Z_2$ on its Lie algebra coincides
with the action of $SO_3$ on $\R^3$, every nonzero
$\pi_1(C)$-invariant vector field $\tilde\xi$ defines an invariant
point on $S^2\cong \p_1$. Therefore $\tilde{\xi} =0$, and $\xi$ is
trivial. \QED

\subsubsection{$c_1(V)=1$ mod 2 }

\begin{defn}
Let $(M, J, \omega)$ be a K\"ahler manifold of complex dimension
$n$. The slope of a holomorphic vector bundle $E$ over $M$ of rank
$r$ is the number $$\mu(E):=\frac{1}{r}\int_M c_1(E)\wedge
\omega^{n-1}.$$
\end{defn}

\begin{defn}
Let $(M, J, \omega)$ be a K\"ahler manifold. A holomorphic vector
bundle is said to be stable if $\mu(F)<\mu(E)$  for any proper
sub-bundle $F\subset E$.
\end{defn}

\begin{defn}
A vector bundle $E$ is said to be polystable  if it decomposes as a
direct sum  of  stable vector bundles with the same slope.
\end{defn}

In \cite{at}, Apostolov and Tonnesen-Friedman have shown that a
complex geometrically ruled surface $M=\p(V)$ over a Riemann surface
$C$, where $V$ is a holomorphic rank 2  vector bundle over $C$,
admits cscK metrics if and only if the bundle $V$ is polystable.

\begin{thm}\label{nontrivialextension}
Let $C$ be a curve of genus at least 1. Let $\mathcal{O}(p)$ be the line bundle over $C$ associated with
the divisor of a point $p\in C$. Then every nontrivial extension $V$ of $\mathcal{O}(p)$ by $\mathcal{O}$ is stable.
\end{thm}
\NI{\it Proof.} (See \cite{f}.) Let $V$ be a nontrivial extension of the form
$$0\to\mathcal{O}\to V\to\mathcal{O}(p)\to 0.$$
(It exists since $H^1(C, \mathcal{O}(-p))=H^0(C, \Omega^1(p))\neq 0$
if $C$ has genus at least 1.) The normalized degree $\mu(V)$ is 1/2.
We need to show that for every subline bundle $L$ of $V$,
$\mu(L)=\deg L$ is less than or equal to zero. If the composite map
$L\to\mathcal{O}(p)$ is zero, then $L$ is contained in
$\mathcal{O}$, and $\deg L\leq 0<\mu(V)$. Otherwise, the map
$L\to\mathcal{O}(p)$ is nonzero. Thus, $L^{-1}\otimes\mathcal{O}(p)$
has a nonzero section, which implies $\deg L\leq 1$. In particular
$\deg L=1$ if and only if $L\cong\mathcal{O}(p)$, and the exact
sequence splits, contrary to the hypothesis. Therefore, we also have
$\deg L\leq 0< \mu(V)$. This completes the proof. \QED

\begin{thm}\label{quasistable2}
Let $C$ be a curve of genus at least 1. Let $\mathcal{O}(p)$ be the
line bundle over $C$ associated with the divisor of a point $p\in
C$. Then the projectivization $S:=\p(V)$ of every nontrivial
extension $V$ of $\mathcal{O}(p)$ by $\mathcal{O}$ admits no
nontrivial holomorphic vector fields which vanish somewhere on $S$.
\end{thm}

\NI{\it Proof.} First consider the case that $C$ has genus at least
2. Let $S:=\p(V)$ be the projectivization of $V$, and $\pi: S\to C$
be the canonical projection. Denote by $Aut(V)$ the automorphism of
$V$ over $C$, by $Aut(C)$ the automorphism group of $C$, and
$Aut_C(S)$ denotes the automorphism group of $S$ over $C$. That is,
$Aut_C(S)=\{\sigma\in Aut(S)|\pi\sigma=\pi\}$. Since $C$ is an
irrational curve, there is an exact sequence
$$ 1 \to Aut_C (S) \to Aut (S) \to Aut (C) .$$  It is well known that every Riemann surface $C$ of genus at least 2
has a discrete automorphism group. Therefore to show that $S$ admits
no nontrivial holomorphic vector fields, it suffices to show the
group $Aut_C(S)$ is discrete. The relation between $Aut(V)$ and
$Aut_C(S)$ can be found in \cite{gr}, which states that one has the
following exact sequence of groups
$$ 1 \to Aut(V)/\Gamma(C,\mathcal{O}^*) \to Aut_C(S)  \to  \triangle\to 1,$$
where $\Gamma(C, \mathcal{O}^*)$ is the group of global holomorphic
sections of the sheaf $\mathcal{O}^*$ over $C$, and
$\triangle:=\{L|L\to C \,\mbox{is a holomorphic line bundle
satisfying}\, V\otimes L=V\}$ is a subgroup of $2$-torsion part of
the Jacobian variety of $C$, hence $\triangle$ is discrete. In the
case that $V$ is indecomposable, Maruyama \cite{m} shows that
$$Aut(V)=\left\{\left( \begin{array}{cc}
\alpha &  s  \\
0 &  \alpha
\end{array} \right)|\, \alpha\in\Gamma(L^*\otimes L)=\C^*, s\in \Gamma(C, (\mbox{det V})^{-1}\otimes L^2)\right\},$$
where $L$ is a maximal sub-bundle of $V$. From the proof of Theorem
\ref{nontrivialextension}, we know the line bundle $(\mbox{det
V})^{-1}\otimes L^2$ has negative degree, and admits no nontrivial
holomorphic sections. This implies $Aut(V)=\C^*=\Gamma(C,
\mathcal{O}^*)$, hence $Aut_C(S)=\triangle$ is discrete.

If the curve $C$ has genus 1, that $Aut_C(S)$ is discrete can be
shown by the same argument. It implies that there exists no vertical
holomorphic vector fields on $S=\p(V)$. Let $\xi$ be a holomorphic
vector field on $S$. Since $S$ is a fibration with every fibre
smooth and compact, $\xi$ projects to a holomorphic vector field
$\bar{\xi}$ on $C$, which can only be parallel. Therefore $\xi$
vanishes nowhere. \QED

\begin{cor} Let $S$ be a ruled surface over a curve $C$ of genus at least 1, then $S$ is deformation equivalent
to a compact complex surface which admits a cscK metric.
\end{cor}

\NI{\it Proof.} First assume  that $S$ is minimal. Since the
deformation equivalence class of $S$ is determined by genus $g(C)$
of $C$ and $c_1(V)$, $S$ is deformation equivalent to either the
surface we constructed in section \ref{representation}, the trivial
bundle $T\times \p_1$ over a torus, or the projectivization of the
nontrivial extension of $\mathcal{O}(p)$ by $\mathcal{O}$. By
Theorem \ref{quasistable} and \ref{nontrivialextension}, each of
them admits cscK metrics. If $S$ is non-minimal, $S$ is deformation
equivalent to the blow-up of the surface we constructed in Section
\ref{representation} or the projectivization of the nontrivial
extension of $\mathcal{O}(p)$ by $\mathcal{O}$. Each of them admits
no nontrivial holomorphic fields by Theorem \ref{quasistable} and
\ref{quasistable2}. Therefore a direct application of Arezzo-Pacard
Theorem $\ref{ap1}$ gives us the conclusion. \QED

The remaining case is the ruled surfaces over a rational curve
$\p_1$. First, assume the surface $S$ over $\p_1$ is minimal. Then
$S$ is either $\p_2, \p_1\times\p_1,$ or
$\p(\mathcal{O}\oplus\mathcal{O}(k))$, where $k\in\N$. The first two
cases admit cscK metrics due to the existence of Fubini-Study metric
on projective space. Denote $\p(\mathcal{O}\oplus\mathcal{O}(k))$ by
$F_k$, then $F_1$ is isomorphic to the one point blow-up of $\p_2$,
and $F_k$ is deformation equivalent to $F_{k'}$ if and only if
$k\equiv k'(mod \,2)$. Since the automorphism group of $\p_2$ is $\p
GL_3(\C)$, if $X$ is the blow-up of $\p_2$ at more than 3 points in
general position, then the automorphism group of $X$ is trivial.

\begin{prop}\label{cscK}
Let $S$ be a compact complex rational surface. Suppose $S$ is not
deformation equivalent to $\p_2\# k\overline{\p}_2$, where $k=1,2$.
Then $S$ is deformation equivalent to a complex surface with cscK
metrics.
\end{prop}
\NI{\it Proof.} The assumption implies that $S$ is either $\p_2$, or
deformation equivalent to $\p_1\times\p_1$, or $\p_2\#
k\overline\p_1$, where $k\geq 3.$ In the first two cases, the
Fubini-Study metric and the product of Fubini-Study metrics provide
a cscK metric. If $S$ is the blow-up of $\p_2$ at $k$ points, where
$3\leq k\leq 8$, Tian and Yau \cite{ty} have shown that $S$ admits a
K\"ahler-Einstein metric. If $X$ is the blow-up of $\p_2$ at more
than $4$ points in the general position, then the automorphism group
$Aut(X)$ is discrete. Therefore a direct application of
Arezzo-Pacard Theorem \ref{ap1} shows that $X=\p_2\#
k\overline{\p}_2$ admits a cscK metric whenever $k\geq 4$. \QED

\subsection{Non-existence case}\

In this subsection, we will show that  $\p_2\# k\overline{\p}_2$,
$k=1,2$, is not deformation equivalent to any complex surface with
cscK metrics. Since the deformation class of $\p_2\#\overline{\p}_2$
consists exactly of Hirzebruch surfaces $F_k=\p(\mathcal
O\oplus\mathcal O(k))$, where $k$ is odd, it suffices to show  the
Lie algebra of holomorphic vector fields $\mathfrak{h}$ of $F_k$ is
not reductive whenever $k>0$.

\begin{defn}
We say a Lie algebra $\mathscr{G}$ is reductive if it is the direct
sum of an abelian and a semisimple Lie algebra.
\end{defn}

Recall the theorem of Lichnerowicz and Matsushima \cite{li, ma}
tells us that a compact K\"{a}hler manifold $(M, J)$ whose identity
component $Aut_0(M,J)$ of the automorphism group is not reductive
does not admit any  cscK metric.

\begin{thm}(Lichnerowicz-Matsushima \cite{li, ma})
Let $(M,J)$ be a compact manifold with cscK metrics. Then the Lie
algebra $\mathfrak{h}(M)$ of holomorphic vector fields decomposes as
a direct sum
$$
\mathfrak{h}(M)=\mathfrak{h}'(M)\oplus\mathfrak{a }(M)
$$
where $\mathfrak{a}(M)$ is the abelian subalgebra of parallel
holomorphic vector fields, and $\mathfrak{h}'(M)$ is the subalgebra
of holomorphic vector fields with zeros. Furthermore,
$\mathfrak{h}'(M)$ is the complexification of the killing fields
with zeros. In particular, $\mathfrak{h}(M)$ is a reductive Lie
algebra.
\end{thm}

Let $E={\mathcal O}\oplus \mathcal O(k), F_k=\p(E)$ be the
projectivization, and $\pi: F_k\rightarrow \p_1$ be the holomorphic
ruling. Let $\mathfrak{ h}(M)$ be the Lie algebra of the holomorphic
vector fields. Then we have the algebra homomorphism
$\phi:\mathfrak{h} \rightarrow \mathfrak{sl}_2(\C)$, and the
following exact sequence:
$$
0 \rightarrow \mathfrak{h}^{\bot} \rightarrow \mathfrak{h}
\xrightarrow{\phi} \mathfrak{sl}_2(\C),
$$
where $\mathfrak{h}^{\bot}$  denotes the algebra of all vertical
holomorphic vector fields. Here we use "vertical" to mean they are
tangent to the fibres.

\begin{lemma}
The algebra homomorphism $\phi:\mathfrak{h} \rightarrow
\mathfrak{sl}_2(\C)$ is surjective.
\end{lemma}

\NI {\it Proof.} Given a holomorphic vector field $\xi\in
\mathfrak{h}(\p_1)$, the generated automorphism group $h_t=exp
(t\xi) \in Aut( \p_1)=\p GL_2(\C)=SL_2(\C)/ \Z_2$ can be lifted to a
family of linear automorphism group $\tilde h_t\in SL_2(\C)$ on
$\C^2$, which fixes the origin $O$, and $\tilde h_0=$ identity.
Using that $\mathcal O(-1)$ is the one point blow-up of $\C^2$ at
the origin and apply Hartog's extension theorem, the derivative
$\tilde \xi=\frac{\partial \tilde h_t}{\partial t}_{|t=0}$
represents a holomorphic vector field on $\mathcal O(-1)$ such that
$p_*\tilde\xi=\xi$, where $p:\mathcal O(-1)\to\p_1$ is the natural
projection. Since $\mathcal O(1)$ is the dual of $\mathcal O(-1)$,
$\mathcal O(k)$ is the tensor product of $k$ copies of $\mathcal
O(1)$, and $\mathcal O\oplus\mathcal O(k)$ is the direct sum of
bundles $\mathcal O$ and $\mathcal O(k)$, we deduce that the family
$h_t$ on $\p_1$ induces a family of automorphisms $\tilde h_t$ on
$E=\mathcal O\oplus\mathcal O(k)$ such that the following diagram
commutes.
$$\begin{CD}
  E  @>\tilde h_t>> E     \\
  @VVV           @VVV \\
   \p_1          @>h_t>>  \p_1 .
\end{CD}$$
Since $\tilde h_t$ restricted to each fibre is linear, it induces a
family of infinitesimal automorphisms on $F_k=\p(E)$, and a
holomorphic vector field $\hat \xi\in\mathfrak{h}(F_k)$ such that
$\pi_* \hat\xi= \xi$. \QED

Recall the exact sequence of groups found by Grothendieck \cite{gr}:
$$1 \to Aut(V)/\Gamma(C,\mathcal{O}^*) \to Aut_C(S)  \to  \triangle\to 1,$$
where $\triangle$ is a discrete subgroup of Picard group. Assume $k>0$.
We have
$$
\begin{array}{ccc}
Aut(E)&=&\Gamma(\p_1, E\otimes E^*) \\
&=& \C\oplus \C\oplus \odot^k V, \\
\end{array}
$$
where $V=\C^2$, and $\odot^k V =\{b_0Z_0^k + b_1
Z_0^{k-1}Z_1+b_2Z_0^{k-2}Z_1^2+...+b_k Z_1^k| b_0, b_1,...,b_k\in
\C\}$ is the set of homogeneous polynomials of degree $k$.  Using
the decomposition $E=\mathcal O\oplus\mathcal O(k)$, we can write
$$
\Gamma(E\otimes E^*)=\left( \begin{array}{cc}
a &  0  \\
b &  c
\end{array} \right)
$$
where $a,c$ are constants, and $b\in \odot^k V$. The group
multiplication is given by
$$
\left( \begin{array}{cc}
a &  0  \\
b &  c
\end{array} \right)
\cdot \left( \begin{array}{cc }
a' &  0  \\
b' &  c'
\end{array} \right)
=\left( \begin{array}{cc }
aa' &  0  \\
ab'+bc' &  cc'
\end{array} \right).
$$
Observe that $aI, a\in\C$ induces the identity map on $\p(E).$
We can then identify the fibre-preserving automorphism group of
$F_k$ with the group
$$
\left\{\left( \begin{array}{cc}
1 &  0  \\
b &  c
\end{array} \right)\mid c\in \C, b\in \odot^kV  \right\}.$$
The Lie algebra $\mathfrak{h}^{\bot}$ of vertical holomorphic vector
fields may be visualized as the set of
$$
\left \{\left( \begin{array}{cc}
0 &  0  \\
\ast &  \times
\end{array} \right)\right \},
$$
where $\times$ stands for an arbitrary complex number, and $\ast$
stands for an arbitrary homogeneous polynomial with degree $k$ in
the variables $z_1,z_2$. A direct computation shows that the center
of $\mathfrak{h}^{\bot}$ is trivial, and $[\mathfrak{h}^{\bot},
\mathfrak{h}^{\bot}] \lvertneqq \mathfrak{h}^{\bot}$. Therefore,
$\mathfrak{h}^{\bot}$ is not reductive.

\begin{thm}
$\mathfrak{h}$ is not reductive. In particular, $F_k$ admits no
cscK metrics whenever $k>0$.
\end{thm}
\NI{\it Proof.} Suppose that $\mathfrak{h}$ is reductive. We claim
that $Z(\mathfrak{h})\subset Z(\mathfrak{h}^{\bot})$, where
$Z(\mathfrak{h}), Z(\mathfrak{h}^{\bot})$ denote the centers of
$\mathfrak{h}$ and $\mathfrak{h}^{\bot}$, respectively. In
particular, $Z(\mathfrak{h})=0$. Recall that we have the exact
sequence of Lie algebras
$$
0 \rightarrow \mathfrak{h}^{\bot} \rightarrow \mathfrak{h}
\xrightarrow{\phi} \mathfrak{sl}_2(\C) \rightarrow 0.
$$
Given $h\in Z(\mathfrak{h})$, $[h, h']=0$ for any $h'\in
\mathfrak{h}$ implies $[\phi(  h), \phi( h')]=0$ in
$\mathfrak{sl}_2(\C)$. Therefore $\phi(h)$ is in the center of
$\mathfrak{sl}_2(\C)$, which is trivial. It follows that $h$ is in
$\mathfrak{h}^{\bot}$, so $h\in Z(\mathfrak{h}^{\bot})$. Since
$\mathfrak{h}$ is reductive, and $\mathfrak{h}^{\bot}$ is a subideal
of $\mathfrak{h}$, there exists a subideal $\tilde {\mathfrak{h}}$
of $\mathfrak {h}$ such that $\mathfrak{h}=\mathfrak{h}^{\bot}\oplus
\tilde{\mathfrak{h}}$. Therefore we have
$$
\begin{array}{ccl}
\mathfrak{h}^{\bot}\oplus\tilde{\mathfrak{h}}&=&\mathfrak{h}\\
&=&[\mathfrak {h}, \mathfrak {h}]\oplus Z(\mathfrak{h})\\
&=&[\mathfrak{h}^{\bot}\oplus\tilde{\mathfrak {h}},\mathfrak
{h}^{\bot}\oplus\tilde{\mathfrak {h}}]\\
&=&[\mathfrak {h}^{\bot},\mathfrak
{h}^{\bot}]\oplus[\tilde{\mathfrak {h}},\tilde{\mathfrak {h}}].
\end{array}
$$
This implies $\mathfrak{h}^{\bot}=[\mathfrak{h}^{\bot},
\mathfrak{h}^{\bot}
]$, which is the contradiction. \QED

Next we consider the deformation class of $\p_2\# 2\overline{\p}_2$. Every
compact complex ruled surface $S$ deformation equivalent to $\p_2
\#2\overline{\p}_2$ can be realized as a one point blow-up of a Hirzebruch
surface $F_k$ at $p$ for some point $p\in F_k$, $k\in \N$. Denote
the surface $F_k\# \overline{\p}_2$ by $S$. Then the Lie algebra
$\mathfrak{h}(S)$ of holomorphic vector fields on $S$ is isomorphic
to the Lie algebra $\mathfrak{h}'(F_k)$ of holomorphic vector fields
vanishing at the point $p$, i.e.
$\mathfrak{h}(S)=\mathfrak{h}'(F_k)=\{X\in\mathfrak{h}(F_k)|
X(p)=0\}$. Consider the algebra homomorphism
$\pi_*:\mathfrak{h}(S)\to \mathfrak{sl}_2(\C)$ composed by
$\mathfrak{h}(S)\to \mathfrak{h}(F_k)$ and
$\mathfrak{h}(F_k)\xrightarrow{\phi} \mathfrak{sl}_2(\C)$. Let
$\mathfrak{g}$ be the image of $\pi _*$. Then  $\mathfrak{g}$ is the
Lie algebra of the subgroup of $SL_2(\C)$ fixing the point $\pi(p)$.
Without loss of generality, we can assume $\pi(p)=[1:0]\in\p_1.$

\begin{thm}
$\mathfrak{h}$(S) is not reductive. In particular, $F_k \# \overline{\p}_2$
admits no cscK metrics if $k\geq 0$.
\end{thm}
\NI{\it Proof.} Suppose $\mathfrak{h}(S)=\mathfrak{h}'(F_k)=\mathfrak{h}'$ is
reductive. Then
$$
\begin{array}{lll}
\mathfrak{g}&=&\pi_*\mathfrak{h}'\\
&=&\pi_*( [\mathfrak{h}',\mathfrak{h}']\oplus Z(\mathfrak{h}') )\\
&=&[\pi_*\mathfrak{h}',\pi_*\mathfrak{h}']+ \pi_*Z(\mathfrak{h}')\\
&\leq& [\mathfrak{g}, \mathfrak{g}]+Z(\mathfrak{g}).
\end{array}
$$
This is absurd since $\mathfrak{g}=\left\{\left( \begin{array}{cc}
a &  b  \\
0 &  -a
\end{array} \right)|a, b \in \C \right\}$ has trivial center, and
$[\mathfrak{g},\mathfrak{g}]= \left\{\left(
\begin{array}{cc}
0 &  \ast  \\
0 &  0
\end{array} \right) | \ast\in \C \right\}$.
\QED

Recently, Arezzo, Pacard, and Singer \cite{ap3} have shown that
there exists extremal metrics on the one or two points blow-up of
$\p_2\# k\overline{\p}_2$. The following corollary follows
immediately from their result and Proposition \ref{cscK}.

\begin{cor}
Let $S$ be a compact complex  surface with even Betti number $b_1$.
Then $S$ is deformation equivalent to a complex surface with
extremal metrics.
\end{cor}

\section{Remarks}
Theorem \ref{general} and \ref{tori} tell us that every compact
complex surface of K\"ahler type and Kodaira dimension $0$ or $2$
carries cscK metrics. This is not true for other cases. As we show
in the nonexistence case, a Hirzebruch surface
$S=\p(\mathcal{O}\oplus \mathcal{O}(k))$ with $k>0$ does not admit
any cscK metrics since its automorphism group is not reductive.
Furthermore, although most of elliptic surfaces and ruled surfaces
with cscK metrics in this article have discrete automorphism group,
and we know there is an $h^{1,1}$-dimensional family of cscK metrics
on the nearby complex surfaces in the deformation class by using
LeBrun and Simanca's theorem \cite{ls}, it is challenging to see,
for a fixed complex structure, which K\"ahler class do these cscK
metrics lie in.   Even if the first Chern class $c_1(X)$ of a
manifold $X$ is negative, there are examples by Ross \cite{ro} where
some K\"ahler classes do not contain any cscK metric.

 \smallskip
\NI Department of Mathematics,\\
South Hall, room 6607,\\
University of California,\\
Santa Barbara, CA 93106, USA . {\small  \newline \noindent\noindent
{\tt yjshu@math.ucsb.edu}}
\end{document}